\documentclass[nonacm]{acmart}
\makeatletter
\def\@ACM@checkaffil{
    \if@ACM@instpresent\else
    \ClassWarningNoLine{\@classname}{No institution present for an affiliation}%
    \fi
    \if@ACM@citypresent\else
    \ClassWarningNoLine{\@classname}{No city present for an affiliation}%
    \fi
    \if@ACM@countrypresent\else
        \ClassWarningNoLine{\@classname}{No country present for an affiliation}%
    \fi
}
\makeatother
\usepackage{graphicx}
\usepackage{multirow}
\usepackage{amsmath}
\usepackage{amsfonts}

\title{Hyb Error: A Hybrid Metric Combining Absolute and Relative Errors}
\author{Peichen Xie}
\affiliation{Microsoft Research}
\authorsaddresses{}

\begin{document}

\begin{abstract}
Suppose $x$ is an approximation of $y$. This paper proposes using $\frac{|x-y|}{1+|y|}$, named Hyb Error, to measure the error. This metric equals half the harmonic mean of absolute error and relative error, effectively combining their advantages while mitigating their limitations. For example, Hyb Error approaches absolute error as $|y|$ approaches 0, thereby avoiding the exaggeration of relative error, and approaches relative error as $|y|$ approaches infinity, thereby avoiding the exaggeration of absolute error. The Hyb Error of $\epsilon$ is equivalent to $|x-y|=\epsilon+\epsilon |y|$, which implies $\mathrm{isclose}(x,y,\epsilon,\epsilon)=\mathrm{True}$, where ``isclose'' is a common floating-point equality check function in numerical libraries. For sequences, this property makes the Maximum Element-wise Hyb Error (MEHE) a pragmatic error metric that reflects the most significant error and equals the decision boundary of the ``isclose'' function.
\end{abstract}

\maketitle

\section{Motivation}

Given an approximate value $x$ and its exact value $y$, the absolute error $|x-y|$ and the relative error $\frac{|x-y|}{|y|}$ are two standard error metrics. However, each metric has limitations that hinder its sole use for error measurement.

Absolute error is not suitable for large values, as it can overstate the importance of error when $|y|$ is large, even if the relative error is small. For example, if $x=1.02$ and $y=1$, the absolute error is 0.02. However, if $x=1020$ and $y=1000$, the absolute error is 20 while the relative error remains 0.02 or 2\%. In this case, the absolute error significantly exaggerates the error, and therefore the relative error is a better metric.

Conversely, relative error is not appropriate for small values, as it exaggerates the error when $|y|$ is small. For example, if $x=1.02$ and $y=1$, the relative error is 2\%. However, if $x=0.021$ and $y=0.001$, the absolute error remains 0.02, but the relative error is 20 or 2000\%, which significantly exaggerates the error. In this case, the absolute error is a better metric.

These examples illustrate the necessity of considering both absolute and relative errors when measuring errors. However, it is more convenient to adopt a single metric. Thus, the question arises: is there a metric that combines the benefits of both absolute and relative errors?

\section{Definition}

In this paper, a hybrid error metric, termed ``Hyb Error'', is introduced. This metric combines the advantages of both absolute and relative errors while mitigating their limitations. With Hyb Error, people can measure errors with a single metric, without overstating errors in any situation.

The design of Hyb Error is derived from the relative error. Because the small value of the denominator in the relative error $\frac{|x-y|}{|y|}$ amplifies the error when $|y|$ is small, it is necessary to enlarge the denominator. Meanwhile, when $|y|$ is large, the magnitude of the denominator should be preserved to retain the benefit of relative error.

To achieve this, the denominator of the relative error is modified by adding one, resulting in
\[\mathrm{HybErr}(x,y)=\frac{|x-y|}{1+|y|}.\]
This paper calls this metric ``Hyb Error''. This design prevents the denominator from being too small, thereby avoiding the exaggeration of errors for small values.

\section{Properties}

The design of adding one is deliberate. Choosing one makes Hyb Error approach absolute error as $|y|$ approaches zero. In this case, one dominates the denominator, resulting in
\[\lim_{|y|\to 0} \mathrm{HybErr}(x,y) = |x-y|. \tag{Property 1}\label{eq:p1}\]

Meanwhile, Hyb Error maintains the magnitude of the denominator when $|y|$ is large, thus avoiding the exaggeration of errors for large values. Because $|y|$ dominates the denominator in this case, Hyb Error approaches relative error as $|y|$ approaches infinity:
\[\lim_{|y|\to +\infty} \mathrm{HybErr}(x,y)=\frac{|x-y|}{|y|}. \tag{Property 2}\label{eq:p2}\]

In addition, Hyb Error is a hybrid metric because it equals half the harmonic mean of absolute error and relative error:
\[\mathrm{HybErr}(x,y)=\frac{1}{2}\times \frac{2}{\frac{1}{|x-y|}+\frac{|y|}{|x-y|}}. \]
This makes Hyb Error a smooth and balanced mixture of absolute and relative errors. Specifically, observing the ratio of Hyb Error to absolute error
\[k_\mathrm{abs}(x,y)=\frac{\mathrm{HybErr}(x,y)}{|x-y|}=\frac{1}{1+|y|}\]
and the ratio of Hyb Error to relative error
\[k_\mathrm{rel}(x,y)=\frac{\mathrm{HybErr}(x,y)}{\frac{|x-y|}{|y|}}=\frac{1}{1+\frac{1}{|y|}},\]
these ratios are smooth and symmetric with respect to the plane $k=0.5$. Formally, 
\[ \forall x \forall y\neq 0 : k_\mathrm{abs}(x,y)+k_\mathrm{rel}(x,y)=1 \land \nabla k_\mathrm{abs} \in \mathbb{R} ^2 \land \nabla k_\mathrm{rel} \in \mathbb{R} ^2. \tag{Property 3}\label{eq:p3}\]

\begin{figure}[h]
    \centering
    \includegraphics[width=0.5\linewidth]{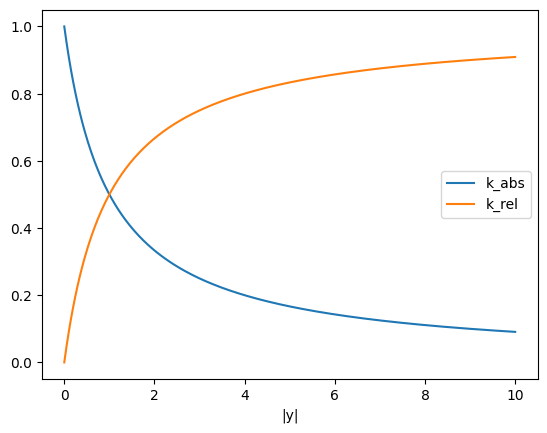}
    \caption{The ratio of Hyb Error to absolute error (blue) and the ratio of Hyb Error to relative error (orange).}
    \label{fig:ratio}
\end{figure}

For Hyb Error, $k_\mathrm{abs}$ and $k_\mathrm{rel}$ are monotonic and smooth functions of $|y|$, as shown in Figure \ref{fig:ratio}. Due to Property 1 and Property 2, it is easy to see that as $|y|$ approaches $0$, $k_\mathrm{rel}$ approaches $0$ and $k_\mathrm{abs}$ approaches $1$. Conversely, as $|y|$ approaches $+\infty$,  $k_\mathrm{rel}$ approaches $1$ and $k_\mathrm{abs}$ approaches $0$.

Next, note that
\[\mathrm{HybErr}(x,y)=\epsilon \Leftrightarrow |x-y|=\epsilon+\epsilon |y|.\]
Therefore, a Hyb Error of $\epsilon$ means that the distance between $x$ and $y$ is $\epsilon+\epsilon |y|$. The distance consists of two parts: the absolute part $\epsilon$ and the relative part $\epsilon |y|$. It denotes the sum of an absolute perturbation of $\epsilon$ and a relative perturbation of $\epsilon$.

The perturbations match the parameter definitions of the ``isclose'' function in popular numerical libraries such as NumPy \cite{harris_array_2020} and PyTorch \cite{paszke_pytorch_2019}. The ``isclose'' function is a common floating-point equality check function defined as
\[\mathrm{isclose}(x,y,a,r)=
\begin{cases}
    \mathrm{True}, & \text{if } |x-y|\leq a+r|y|\\
    \mathrm{False}, & \text{otherwise}
\end{cases},\]
where the parameter $a$ represents absolute tolerance and the parameter $r$ represents relative tolerance. Therefore, the Hyb Error of $\epsilon$ implies that $x$ is close to $y$ at an absolute tolerance level of $\epsilon$ and a relative tolerance level of $\epsilon$. Formally,
\[\mathrm{HybErr}(x,y) = \epsilon \Rightarrow \mathrm{isclose}(x,y,\epsilon,\epsilon) = \mathrm{True}, \]
and 
\[\mathrm{HybErr}(x,y) = \epsilon \Rightarrow \forall \epsilon' \in [0,\epsilon) : \mathrm{isclose}(x,y,\epsilon',\epsilon') = \mathrm{False}. \]
In other words,
\[\mathrm{HybErr}(x,y) = \min \{\epsilon: \mathrm{isclose}(x,y,\epsilon,\epsilon) = \mathrm{True}\} \tag{Property 4}\label{eq:p4}\]

Lastly, unlike relative error, which is undefined when $y=0$, Hyb Error is well defined for all $x\in \mathbb{R}$ and $y\in \mathbb{R}$. Specifically, 
\[\mathrm{HybErr}(x,0)=|x| \tag{Property 5}\label{eq:p5}\]

Table \ref{tab:properties} summarizes the properties.

\begin{table}[h]
    \centering
    \begin{tabular}{|c|c|} \hline 
         & Description\\ \hline 
 Property 1&The error approaches absolute error as $|y|\to 0$.\\ \hline 
         Property 2& The error approaches relative error as $|y|\to +\infty$.\\ \hline 
         Property 3& The error is a smooth and balanced mixture of absolute and relative errors.\\ \hline 
         Property 4& The error is the decision boundary of $\epsilon \mapsto \mathrm{isclose}(x,y,\epsilon,\epsilon)$.\\ \hline 
         Property 5& The error is well defined when $y=0$.\\ \hline
    \end{tabular}
    \caption{Properties of Hyb Error.}
    \label{tab:properties}
\end{table}

\section{Discussion}

This section discusses alternative designs and shows that Hyb Error is better than them.

\subsection{Adding other numbers}

Hyb Error can be treated as a smoothed version of relative error with the smooth factor of one added to the denominator. Naturally, people may choose another smooth factors, resulting in the first alternative design
\[\frac{|x-y|}{t+|y|}\]
where the smooth factor $t\neq 1$.

However, this error metric no longer has Property 1, because
$\lim_{|y|\to 0} \frac{|x-y|}{t+|y|} = \frac{|x-y|}{t}\neq |x-y|$.

\begin{figure}[h]
    \centering
    \includegraphics[width=0.5\linewidth]{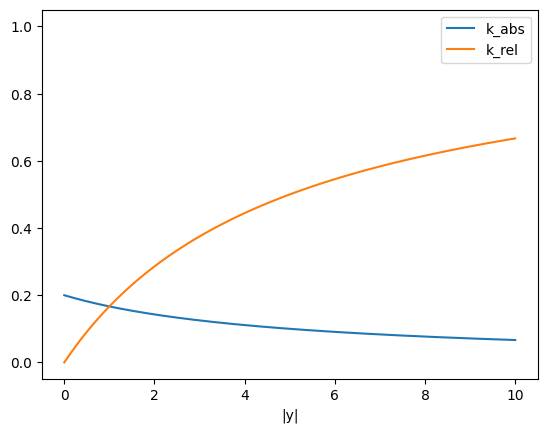}
    \caption{The ratio of $\frac{|x-y|}{5+|y|}$ to absolute error (blue) and to relative error (orange).}
    \label{fig:ratio1}
\end{figure}

It does not have Property 3 either, because the ratio $k_\mathrm{abs}=\frac{1}{t+|y|}$ and the ratio $k_\mathrm{rel}=\frac{1}{1+t/|y|}$ are now asymmetric, as shown in Figure \ref{fig:ratio1}. 

It does not have Property 4 either, but has a variant version: $\frac{|x-y|}{t+|y|} = \min \{\epsilon: \mathrm{isclose}(x,y,t\epsilon,\epsilon) = \mathrm{True}\}$.

\subsection{Using maximum}

\begin{figure}[h]
    \centering
    \includegraphics[width=0.5\linewidth]{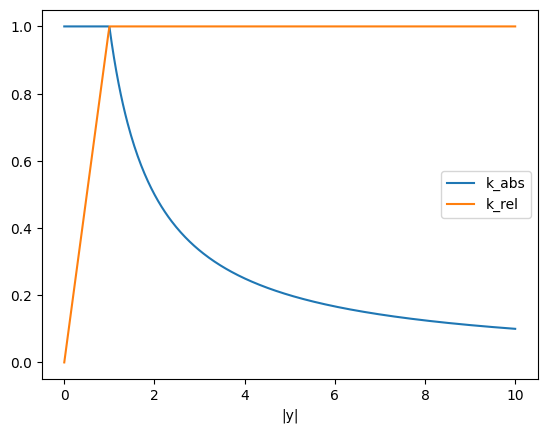}
    \caption{The ratio of $\frac{|x-y|}{\max(1,|y|)}$ to absolute error (blue) and to relative error (orange).}
    \label{fig:ratio2}
\end{figure}

Change the addition operation in Hyb Error to the maximum operation, the second alternative design
\[\frac{|x-y|}{\max(1,|y|)}\]
has both Property 1 and Property 2. However, it does not have Property 3, because now $k_\mathrm{abs} = \frac{1}{\max(1,|y|)}$ and $k_\mathrm{rel} = \frac{|y|}{\max(1,|y|)}$ are non-smooth and asymmetric functions, as shown in Figure \ref{fig:ratio2}. 

It does not have Property 4 either, but it has a variant version $\frac{|x-y|}{\max(1,|y|)} = \min \{\epsilon: \mathrm{isclose_{max}}(x,y,\epsilon,\epsilon) = \mathrm{True}\}$ where
\[\mathrm{isclose_{max}}(x,y,a,r)=
\begin{cases}
    \mathrm{True}, & \text{if } |x-y|\leq \max(a,r|y|)\\
    \mathrm{False}, & \text{otherwise}
\end{cases}\]
is a variant of the isclose function.

In summary, Table \ref{tab:comparison} compares these alternatives with Hyb Error. Hyb Error retains all five desirable properties, whereas the alternative designs lack certain properties.

\begin{table}[h]
    \centering
    \begin{tabular}{|c|c|c|c|c|c|} \hline 
         &  Property 1&  Property 2&  Property 3&  Property 4& Property 5\\ \hline 
         Hyb Error&  yes&  yes&  yes&  yes& yes\\ \hline 
         Alternative 1 ($\frac{|x-y|}{t+|y|}$)&  no&  yes&  no&  no& yes\\ \hline 
         Alternative 2 ($\frac{|x-y|}{\max(1,|y|)}$)&  yes&  yes&  no&  no& yes\\ \hline
    \end{tabular}
    \caption{Comparison between error metrics.}
    \label{tab:comparison}
\end{table}

\section{Errors of vectors and sequences}

\subsection{Hyb Error of vectors}

For a vector $\mathbf{y}=(y_1,y_2,...,y_n)$ and its approximation $\mathbf{x}=(x_1,x_2,...,x_n)$, absolute error $||\mathbf{x}-\mathbf{y}||$ and relative error $\frac{||\mathbf{x}-\mathbf{y}||}{||\mathbf{y}||}$ are standard error metrics, where $||\cdot||$ donates any vector norm. Similarly, the Hyb Error of vectors is defined as
\[\mathrm{HybErr}(\mathbf{x},\mathbf{y})=\frac{||\mathbf{x}-\mathbf{y}||}{1+||\mathbf{y}||}.\]

\subsection{Mean Hyb Error (MHE)}

A more general question is to measure the error between paired sequences. For a sequence  $Y=(y_1,y_2,...,y_n)$ and its approximation $X=(x_1,x_2,...,x_n)$, in addition to treating them as vectors, it is common to calculate element-wise errors and obtain a sequence of errors, and then calculate the magnitude of the errors.

Specifically, mean absolute error $\mathrm{MAE}=\sum_{i=1}^{n}|x_i-y_i|$ and mean relative error $\mathrm{MRE}=\sum_{i=1}^{n}\frac{|x_i-y_i|}{|y_i|}$ (also known as mean absolute percentage error or MAPE in statistics) are standard metrics representing the average magnitude of errors. This can be analogized to the Mean Hyb Error:
\[\mathrm{MHE}=\sum_{i=1}^{n}\frac{|x_i-y_i|}{1+|y_i|}.\]

\subsection{Maximum Element-wise Hyb Error (MEHE)}

A more interesting question is how to measure the maximum error between two sequences. The standard maximum element-wise absolute error $\mathrm{MEAE}=\max_i{|x_i-y_i|}$ and maximum element-wise relative error $\mathrm{MERE}=\max_i \frac{|x_i-y_i|}{|y_i|}$ do not work well due to the error exaggeration issues for very large and very small values.

Instead, the Maximum Element-wise Hyb Error (MEHE), defined as
\[\mathrm{MEHE}=\max_i \frac{|x_i-y_i|}{1+|y_i|},\]
works very well because Hyb Error solves the issues of error exaggeration for both small and large values. This will be exemplified in Section \ref{sec:example}.

Deduced from Property 4 for Hyb Error, MEHE must equal the decision boundary of $\epsilon \mapsto \mathrm{isclose}(X,Y,\epsilon,\epsilon)$. In other words,
\[\mathrm{MEHE}(X,Y)=\min \{\epsilon : \mathrm{isclose}(X,Y,\epsilon,\epsilon)=\mathrm{True}\}. \tag{Property 4'}\label{eq:p4'}\]

Therefore, MEHE reflects the maximum error between two sequences, implying that $X$ is within the neighborhood of $Y$ for an absolute perturbation of MEHE and a relative perturbation of MEHE, as $\forall i : \mathrm{isclose}(x_i,y_i,\mathrm{MEHE},\mathrm{MEHE})=\mathrm{True}$.

 \section{Example}
 \label{sec:example}

This section demonstrates that standard error metrics do not work well while Hyb Error works well consistently. Suppose the sequence $Y$ equals (0.001, 0.001, 0.001, 1, 1, 1, 1000, 1000, 1000) and its approximation $X$ equals (0.021, 0.00102, 0.2102, 1.02, 1.02, 1.04, 1000.02, 1020, 1020.02). In fact, $x_1, x_4, x_7$ are obtained by adding an absolute perturbation of 0.02 to $y_1, y_4, y_7$; $x_2, x_5, x_8$ are obtained by adding a relative perturbation of 0.02 to $y_2, y_5, y_8$; $x_3,x_6,x_9$ are obtained by adding both an absolute perturbation of 0.02 and a relative perturbation of 0.02 to $y_3, y_6, y_9$.

The element-wise errors in different metrics are calculated and shown in Table \ref{tab:example}. When $i = 1, 2, 3$, the absolute values of $y_i$ are small, leading to exaggeration of the relative errors. Meanwhile, the Hyb Errors, which approximate to the absolute errors, represent the magnitude of errors well. 

When $i=7, 8, 9$, the absolute values of $y_i$ are large, leading to exaggeration of the absolute errors. Meanwhile, the Hyb Errors, which approximate to the relative errors, represent the magnitude of errors well.

In these examples, the Hyb Errors effectively avoid the exaggerated relative errors for small values and the exaggerated absolute errors for large values, capturing the significance of the error indicated by the other metric. Indeed, Hyb Error is no greater than the minimum of absolute error and relative error, because $k_\mathrm{abs} \le 1$ and $k_\mathrm{rel} < 1$.

The bottom of Table \ref{tab:example} lists nine different error measurements for $X$ against $Y$: three mean errors, three maximum element-wise errors, and three errors of vectors (using the 1-norm). Although $X$ is obtained by adding perturbations of only 0.02, the mean absolute error (MAE), mean relative error (MRE), maximum element-wise absolute error (MEAE), maximum element-wise relative error (MERE), and absolute error of vectors are all at high magnitudes, because they are affected by exaggerated errors of some elements. Only the relative error of vectors and all three metrics of Hyb Error are at a magnitude about 0.02, which correctly reflects the overall significance of the error. Particularly, the Maximum Element-wise Hyb Error (MEHE) captures the level of perturbation exactly because of Property 4'.

\begin{table}[h]
\centering
\begin{tabular}{@{}llllll@{}}
\toprule
$i$ & $x_i$ & $y_i$ & Absolute error & Relative error & Hyb Error \\ \midrule
1 & 0.021 & 0.001 & 0.02 & 20 & 0.01998 \\
2 & 0.00102 & 0.001 & 0.00002 & 0.02 & 0.00001998 \\
3 & 0.02102 & 0.001 & 0.02002 & 20.02 & 0.02 \\
4 & 1.02 & 1 & 0.02 & 0.02 & 0.01 \\
5 & 1.02 & 1 & 0.02 & 0.02 & 0.01 \\
6 & 1.04 & 1 & 0.04 & 0.04 & 0.02 \\
7 & 1000.02 & 1000 & 0.2 & 0.0002 & 0.0001998 \\
8 & 1020 & 1000 & 20 & 0.02 & 0.01998 \\
9 & 1020.02 & 1000 & 20.02 & 0.02002 & 0.02 \\ \midrule
 & Mean &  & 4.482227 & 4.462247 & 0.013353 \\
 & Max &  & 20.02 & 20.02 & 0.02 \\
 & Vector &  & 40.34004 & 0.013433 & 0.013429 \\ \bottomrule
\end{tabular}
\caption{Measuring errors with different metrics.}
\label{tab:example}
\end{table}

\section{Conclusion}

This paper presents Hyb Error, an error metric with five desirable properties. It can be used as a single metric to measure errors substituting absolute error and relative error. Moreover, Maximum Element-wise Hyb Error (MEHE) is a useful metric reflecting the maximum error between two sequences.

\bibliographystyle{ACM-Reference-Format}
\bibliography{reference}

\end{document}